\documentclass[12pt]{svmult}
\usepackage{amsmath,amssymb,bbm,graphicx,natbib}

\smartqed

\title*{Best attainable rates of convergence for the estimation of the
  memory parameter}
\titlerunning{Best attainable rates of convergence}

\author{Philippe Soulier} \institute{Universit\'e Paris
  Ouest-Nanterre, 200 avenue de la R\'epublique, 92000 Nanterre cedex,
  France \texttt{philippe.soulier@u-paris10.fr} }

\begin{document}

\maketitle

\abstract {The purpose of this note is to prove a lower bound for the
  estimation of the memory parameter of a stationary long memory
  process. The memory parameter is defined here as the index of
  regular variation of the spectral density at 0. The rates of
  convergence obtained in the literature assume second order regular
  variation of the spectral density at zero. In this note, we do not
  make this assumption, and show that the rates of convergence in this
  case can be extremely slow. We prove that the log-periodogram
  regression (GPH) estimator achieves the optimal rate of convergence
  for Gaussian long memory processes}

\section{Introduction}
Let $\{X_n\}$ be a weakly stationary process with autocovariance
function $\gamma$. Its spectral density $f$, when it exists, is an
even nonnegative measurable function such that
\begin{gather*}
  \gamma(k) = \int_{-\pi}^\pi f(x) \mathrm e^{\mathrm i k x} \,
  \mathrm d x \; .
\end{gather*}
Long memory of the weakly stationary process $\{X_t\}$ means at least
that the autocovariance function is not absolutely summable. This
definition is too weak to be useful. It can be strengthen in several
ways. We will assume here that the spectral density is regularly
varying at zero with index $-\alpha \in(-1,1)$, i.e. it can be
expressed for $x\geq0$  as
\begin{gather*}
  f(x) = x^{-\alpha} L(x) \;, 
\end{gather*}
where the function $L$ is slowly varying at zero, which means that for
all $t>0$,
\[
\lim_{x\to 0} \frac{L(tx)}{L(x)} = 1 \;.
\]
Then the autocovariance function is regularly varying at infinity with
index $\alpha-1$ and non absolutely summable for $\alpha>0$. The main
statistical problem for long memory processes is the estimation of the
memory parameter $\alpha$. This problem has been exhaustively studied
for the most familiar long memory models: the fractional Gaussian
noise and the ARFIMA($p,d,q$) process. The most popular estimators are
the GPH estimator and the GSE estimator, first introduced respectively
by \cite{geweke:porter-hudak:1983} and \cite{kunsch:1987}.  Rigorous
theoretical results for these estimators were obtained by
\cite{robinson:1995g,robinson:1995l}, under an assumption of second
order regular variation at 0, which roughly means that there exists
$C,\rho>0$ such that
\begin{gather*}
  f(x) = C x^{-\alpha}\{1+O(x^\rho)\} \; .
\end{gather*}
Under this assumption, \cite{giraitis:robinson:samarov:1997} proved
that the optimal rate of convergence of an estimator based on a sample
of size $n$ is of order $n^{2\rho/(2\rho+1)}$. 

The methodology to prove these results is inspired from similar
results in tail index estimation. If $F$ is a probability distribution
function on $(-\infty,\infty)$ which is second order regularly varying
at infinity, i.e. such that
\begin{gather*}
  \bar F(x) = C x^{-\alpha} \{1+O(x^{-\alpha\rho})\} 
\end{gather*}
as $x\to\infty$, then \cite{hall:welsh:1984} proved that the best
attainable rate of convergence of an estimator of the tail index
$\alpha$ based on $n$ i.i.d. observations drawn from the distribution
$F$ is of order $n^{2\rho/(2\rho+1)}$.  In this context,
\cite{drees:1998} first considered the case where the survival
function $\bar F$ is regularly varying at infinity, but not
necessarily second order regularly varying. He introduced very general
classes of slowly varying functions for which optimal rates of
convergence of estimators of the tail index can be computed. The main
finding was that the rate of convergence can be extremely slow in such
a case.

In the literature on estimating the memory parameter, the possibility
that the spectral density is not second order regularly varying has
not yet been considered. Since this has severe consequences on the
estimations procedures, it seems that this problem should be
investigated. In this note, we parallel the methodology developped by
\cite{drees:1998} to deal with such regularly varying functions. Not
surprisingly, we find the same result, which show that the absence of
second order regular variation of the spectral density has the same
drastic consequences.

The rest of the paper is organised as follows. In
Section~\ref{sec:lowerbounds}, we define the classes of slowly varying
functions that will be considered and prove a lower bound for the rate
of convergence of the memory parameter. This rate is proved to be
optimal in Section~\ref{sec:upperbound}. An illustration of the
practical difficulty to choose the bandwidth parameter is given in
Section~\ref{sec:bandwidth}. Technical lemmas are deferred to
Section~\ref{sec:lemmes}.

\section{Lower bound}\label{sec:lowerbounds}
In order to derive precise rates of convergence, it is necessary to
restrict attention to the class of slowly varying functions referred
to by \cite{drees:1998} as {\em normalised}.  This class is also
referred to as the Zygmund class. Cf.
\cite[Section~1.5.3]{bingham:goldie:teugels:1989}

\begin{definition}
  Let $\eta^*$ be a non decreasing function on $[0,\pi]$, regularly
  varying at zero with index $\rho\geq 0$ and such that
  $\lim_{x\to0}\eta^*(x)=0$.  Let $SV(\eta^*)$ be the class of even
  measurable functions $L$ defined on $[-\pi,\pi]$ which can be for
  expressed $x\geq0$ as
\[
L(x) = L(\pi) \exp\left \{-\int_x^\pi \frac{\eta(s)}s \,ds \right\} \;
,
\]
for some measurable function $\eta$ such that $|\eta| \leq \eta^*$.
\end{definition}
This representation implies that $L$ has locally bounded variations
and $\eta(s) = sL'(s)/L(s)$. Usual slowly varying functions, such as
power of logarithms, iterated logarithms are included in this class,
and it easy to find the corresponding $\eta$ function. Examples are
given below. We can now state our main result.

\begin{theorem}\label{theo:lowerbound}
  Let $\eta^*$ be a non decreasing function on $[0,\pi]$, regularly
  varying at~0 with index $\rho\geq 0$ and such that
  $\lim_{x\to0}\eta^*(x)=0$.  Let $t_n$ be a sequence satisfying
\begin{gather} \label{eq:condition-tn}
\lim_{n\to\infty} \eta^*(t_n) (nt_n)^{1/2} = 1 \; .
\end{gather}
Then, if $\rho>0$, 
\begin{gather}\label{eq:uniformerhopositive}
  \liminf_{n\to\infty} \inf_{\hat \alpha_n} \sup_{L\in SV(\eta^*)}
  \sup_{\alpha\in(-1,1)} \mathbb{E}_{\alpha,L}[\eta^*(t_n)^{-1}|\hat
  \alpha_n - \alpha|] >0 \; ,
\end{gather}
and if $\rho=0$
\begin{gather} \label{eq:uniformeslow}
  \liminf_{n\to\infty} \inf_{\hat \alpha_n} \sup_{L\in SV(\eta^*)}
  \sup_{\alpha\in(-1,1)} \mathbb{E}_{\alpha,L} [\eta^*(t_n)^{-1} |\hat
  \alpha_n - \alpha| ] \geq 1 \; ,
\end{gather}
where $\mathbb{P}_{\alpha,L}$ denotes the distribution of any second order
stationary process with spectral density $x^{-\alpha}L(x)$ and the infimum
$\inf_{\hat\alpha_n }$ is taken on all estimators of $ \alpha$ based on $n$
observations of the process.
\end{theorem}

\begin{example} \label{xmpl:grs}
  Define $\eta^*(s) = C s^\beta$ for some $\beta>0$ and $C>0$. Then
  any function $L\in SV(\eta^*)$ satisfies $L(x) = L(0) +
  O(x^\beta)$, and we recover the case considered
  by~\cite{giraitis:robinson:samarov:1997}. The lower bound for the
  rate of convergence is $n^{\beta/(2\beta+1)}$.
\end{example}

\begin{example} \label{xmpl:logrho}
  For $\rho>0$, define $\eta^*(s) = \rho/\log(1/s)$, then
  $$
  \exp\left\{\int_x^{1/e}\frac{\eta^*(s)}s \,ds\right\} =
  \exp\left\{\rho \log\log(1/x) \right\} = \log^\rho(1/x).
  $$
  A suitable sequence $t_n$ must satisfy $\rho^2/\log^2(t_n)
  \approx nt_n$.  One can for instance choose $t_n =
  \log^2(n)/(n\rho^2)$, which yields $\eta^*(t_n) =
  \rho/\log(n)\{1+o(1)\}$.  Note that $\eta(s) = \rho/\log(s)$ belongs
  to $SV(\eta^*)$, and the corresponding slowly varying function is
  $\log^{-\rho}(1/x)$.  Hence, the rate of convergence is not affected
  by the fact that the slowly varying function vanishes or is infinite
  at~0.
\end{example}
\begin{example}
  The function $L(x) = \log\log(1/x)$ is in the class $SV(\eta^*)$
  with $\eta^*(x) = \{\log(1/x)\log\log(1/x)\}^{-1}$.  In that case,
  the optimal rate of convergence is $\log(n)\log\log(n)$. Even though
  the slowly varying function affecting the spectral density at zero
  diverges very weakly, the rate of convergence of any estimator of
  the memory parameter is dramatically slow.
\end{example}
\noindent{\bf Proof of Theorem \ref{theo:lowerbound}}
  Let $\ell>0$, $t_n$ be a sequence that satisfies the assumption of
  Theorem \ref{theo:lowerbound}, and define $\alpha_n = \eta^*(\ell
  t_n)$ and
\begin{gather*}
  \eta_n(s) = \left\{\begin{array}{ll}
      0 & \mbox{ if $0 \leq s \leq  \ell t_n$, }\\
      \alpha_n & \mbox{ if $ \ell t_n < s \leq \pi$, }
\end{array} \right. \\
L_n(x) = \pi^{\alpha_n} \exp \left\{-\int_x^\pi \eta_n(s) \,ds
\right\}.
\end{gather*}
Since $\eta^*$ is assumed non decreasing, it is clear that $L_n \in
SV(\eta^*)$.  Define now $f_n^-(x) = x^{-\alpha_n}L_n(x)$ and $f_n^+ =
(f_n^-)^{-1}$.  $f_n^-$ can be written as
\[
f_n^-(x) = \left\{ \begin{array}{ll}
    (\ell t_n/x)^{\alpha_n} & \mbox{ if $0 < x \leq \ell t_n$, }\\
    1 & \mbox{ if $\ell t_n < x \leq \pi$. }
\end{array}\right.
\]
 Straighforward computations yield
 \begin{gather}\label{eq:camarche}
 \int_0^\pi \{f_n^-(x) - f_n^+(x)\}^2 \, \mathrm d x 
 = 8 \ell t_n \alpha_n^2( 1 + O(\alpha_n^2)) = 8\ell n^{-1}(1+ o(1)).
 \end{gather}
The last equality holds by definition of the sequence $t_n$.
Let $\mathbb{P}_n^-$ and $\mathbb{P}_n^+$ denote the distribution of 
a $n$-sample of a  stationary
Gaussian processes with spectral densities $f_n^-$ et $f_n^+$
respectively, $\mathbb{E}_n^-$ and $\mathbb{E}_n^+$ the expectation with respect to
these probabilities, $\frac{d\mathbb{P}_n^+}{d\mathbb{P}_n^-}$ the likelihood ratio and
$A_n = \{ \frac{d\mathbb{P}_n^+}{d\mathbb{P}_n^-} \geq \tau\}$ for some real
$\tau\in(0,1)$.  Then, for any estimator $\hat \alpha_n$, based on the
observation $(X_1,\dots,X_n)$, 
\begin{align*}
  \sup_{\alpha,L}  \mathbb{E}_{\alpha,L} [|\hat\alpha_n-\alpha|] & \geq
  \frac12 \left( \mathbb{E}_n^+[|\hat\alpha_n - \alpha_n|]+
    \mathbb{E}_n^-[|\hat\alpha_n + \alpha_n|]\right) \\
  & \geq \frac12 \mathbb{E}_n^- \left[\mathbbm{1}_{A_n}|\hat\alpha_n + \alpha_n| +
    \frac{d\mathbb{P}_n^+}{d\mathbb{P}_n^-} \mathbbm{1}_{A_n}|\hat\alpha_n - \alpha_n| \right] \\
  & \geq \frac12 \mathbb{E}_n^- \left[\{|\hat\alpha_n + \alpha_n| + \tau
    |\hat\alpha_n - \alpha_n|\}\mathbbm{1}_{A_n} \right] \geq \tau \alpha_n
  \mathbb{P}_n^-(A_n).
\end{align*}
Denote $\epsilon=\log(1/\tau)$ and $\Lambda_n =
\log(d\mathbb{P}_n^+/d\mathbb{P}_n^-)$.  Then $\mathbb{P}_n^-(A_n) = 1 - \mathbb{P}_n^-(\Lambda_n
\leq - \epsilon)$. Applying~(\ref{eq:camarche})
and~\cite[Lemma~2]{giraitis:robinson:samarov:1997}, we obtain that
there exist constants $C_1$ and $C_2$ such that
\begin{gather*}
  \mathbb{E}_{n}^-[\Lambda_n] \leq C_1 \ell \; , \ \ 
  \mathbb{E}_{n}^-[(\Lambda_n-m_n)^2] \leq C_2 \ell \; .
\end{gather*}
This yields, for any $\eta>0$ and small enough $\ell$, 
\begin{gather*}
  \mathbb{P}_n^-(A_n) \geq 1 - \epsilon^{-2} \mathbb{E}[\Lambda_n^2] \geq 1 - C
  l\epsilon^{-2} \geq 1 - \eta \; .
\end{gather*}
Thus, for any $\eta, \tau\in(0,1)$, and sufficiently small $\ell$, we
have
\begin{multline*}
  \liminf_{n\to\infty} \inf_{L\in SV(\eta^*)} \inf_{\alpha\in(-1,1)}
  \mathbb{E}_{\alpha,L}[\eta^*(t_n)^{-1}|\hat \alpha_n - \alpha|] \\
  \geq \tau(1-\eta) \lim_{n\to\infty} \frac{\eta^*(\ell
    t_n)}{\eta^*(t_n)} = \tau (1-\eta) \ell^\rho \; .
\end{multline*}
This proves~(\ref{eq:uniformerhopositive})
and~(\ref{eq:uniformeslow}). \qed

\section{Upper bound} \label{sec:upperbound}
In the case $\eta^*(x) = Cx^\rho$ with $\rho>0$,
\cite{giraitis:robinson:samarov:1997} have shown that the lower
bound~(\ref{eq:uniformerhopositive}) is attainable. The extension of
their result to the case where $\eta^*$ is regularly varying with
index $\rho>0$ (for example to functions of the type $x^\beta
\log(x)$) is straightforward.  We will restrict our study to the case
$\rho=0$, and will show that the lower bound~(\ref{eq:uniformeslow})
is asymptotically sharp, i.e.  there exist estimators that are rate
optimal up to the exact constant.

Define the discrete Fourier transform and the periodogram ordinates of
a process $X$ based on a sample $X_1,\ldots,X_n$, evaluated at the
Fourier frequencies $x_j= 2j\pi/n$, $j= 1, \ldots, n$, respectively by
\[
d_{X,j} = (2\pi n)^{-1/2} \sum_{t=1}^n X_t \mathrm{e}^{-\mathrm{i}tx_j}, \ 
\mbox{ and } \ I_{X,j} = |d_{X,j}|^2.
\]
The frequency domain estimates of the memory parameter $\alpha$ are
based on the following heuristic approximation: the renormalised
periodogram ordinate $I_{X,j}/f(x_j)$, $1\leq j \leq n/2$ are
approximately i.i.d. standard exponential random variables.  Although
this is not true, the methods and conclusion drawn from these
heuristics can be rigourously justified.  In particular, the Geweke
and Porter-Hudak (GPH) and Gaussian semiparametric estimator have been
respectively proposed by \cite{geweke:porter-hudak:1983} and
\cite{kunsch:1987}, and a theory for them was obtained
by~\cite{robinson:1995l, robinson:1995g} in the case where the
spectral density is second order regularly varying at 0.

The GPH estimator is based on an ordinary least square regression of
$\log(I_{X,k})$ on $\log(k)$ for $k=1,\ldots,m$, where $m$ is a
bandwith parameter:
\[
(\hat \alpha(m),\hat C) = \arg\min_{\alpha,C} \sum_{k=1}^m
\left\{\log(I_{X,k})- C + \alpha \log(k) \right\}^2.
\]
The GPH estimator has an explicit expression as a weighted 
sum of log-periodogram ordinates:
\[
\hat \alpha(m) = -s_m^{-2}\sum_{k=1}^m \nu_{m,k} \log(I_{X,k}),
\]
with $\nu_{m,k} = \log(k) - m^{-1}\sum_{j=1}^m \log(j)$ and $s_m^2 =
\sum_{k=1}^m \nu_{m,k}^2 = m\{1+o(1)\}$.

\begin{theorem}\label{theo:gph}
  Let $\eta^*$ be a non decreasing slowly varying function such that
  $\lim_{x\to0}\eta^*(x)=0$. Let $\mathbb{E}_{\alpha,L}$ denote the expectation
  with respect to the distribution of a {\em Gaussian} process with spectral
  density $x^{-\alpha}L(x)$. Let $t_n$ be a sequence that
  satisfies~(\ref{eq:condition-tn}) and let $m$ be a non decreasing sequence of
  integers such that
\begin{gather}\label{eq:condmgph}
\lim_{n\to\infty} m^{1/2} \eta^*(t_n) = \infty \mbox{ and }
\lim_{n\to\infty} \frac{ \eta^*(t_n)}{\eta^*(m/n)} =1 \; .
\end{gather}
Assume also that the sequence $m$ can be chosen in such a way that
\begin{gather} \label{eq:condition-m-technique}
  \lim_{n\to\infty} \frac{\log(m) \int_{m/n}^\pi s^{-1}\eta^*(s) \,
    \mathrm d s}{m\eta^*(m/n)} = 0 \; .
\end{gather}
Then, for any  $\delta\in(0,1)$, 
\begin{gather} \label{eq:bornesup}
\limsup_{n\to\infty}\sup_{|\alpha|\leq \delta}\sup_{L\in
  SV(\eta^*)} \eta^*(t_n)^{-2}
\mathbb{E}_{\alpha,L}[(\hat\alpha(m)-\alpha)^2 ] \leq 1.
\end{gather}
\end{theorem}
\begin{remark}
  Since $\eta^*$ is slowly varying, it is always possible to choose
  the sequence $m$ in such a way that~(\ref{eq:condmgph}) holds.
  Condition~(\ref{eq:condition-m-technique}) ensures that the bias of
  the estimator is of the right order. It is very easily checked and
  holds for all the examples of usual slowly varying function $\eta*$,
  but we have not been able to prove that it always holds.

\end{remark}

Since the quadratic risk is greater than the $L^1$ risk, we obtain the
following corollary. 

\begin{corollary}    
  Let $\delta\in(0,1)$ and $\eta^*$ be a non decreasing slowly varying
  function such that $\lim_{x\to0}\eta^*(x)=0$ and such that it is
  possible to choose a sequence $m$ that
  satisfies~(\ref{eq:condition-m-technique}). Then, for $t_n$ as
  in~(\ref{eq:condition-tn}),
\begin{gather} \label{eq:sharp}
  \liminf_{n\to\infty} \inf_{\hat \alpha_n} \sup_{L\in SV(\eta^*)}
  \sup_{\alpha\in(-\delta,\delta)} \mathbb{E}_{\alpha,L}[\eta^*(t_n)^{-1}|\hat
  \alpha_n - \alpha|] = 1 \; .
\end{gather}
\end{corollary}
\begin{remark} 
  This corollary means that the GPH estimator achieves the optimal rate of
  convergence, up to the exact constant over the class $SV(\eta^*)$ when
  $\eta^*$ is slowly varying.  This implies in particular that, contrarily to
  the second order regularly varying case, there is no loss of efficiency of
  the GPH estimator with respect to the GSE.  This happens because in the
  slowly varying case, the bias term dominates the stochastic term if the
  bandwidth parameter $m$ satisfies \eqref{eq:condmgph}. This result is not
  completely devoid of practical importance, since when the rate of convergence
  of an estimator is logarithmic in the number of observations, constants do
  matter.
\end{remark}

\begin{example}[Example \ref{xmpl:logrho} continued] \label{xmpl:convproba} If
  $L(x) = \log^\rho(1/x)\tilde L(x)$, where $\tilde L\in SV(Cx^\beta)$ for some
  $\rho>$, $\beta>0$ and $C>0$, then $\sum_{k=1}^m \nu_{m,k} \log(L(x_k))$
  $\sim \rho m \log^{-1}(x_m)$.  Choosing $m=\log^{1+\delta}(n)$
  yields~(\ref{eq:condmgph}),~(\ref{eq:condition-m-technique}) and
  $\log(n)(\hat\alpha(m)-\alpha)$ converges in probability to $\rho$.

\end{example}

\noindent{\em Proof of Theorem \ref{theo:gph}.} \ 
Define $\mathcal{E}_k = \log\{x_k^\alpha I_k/L(x_k)\}$.  The deviation
of the GPH estimator can be split into a stochastic term and a bias
term:
\begin{gather} \label{eq:devgph}
\hat \alpha(m)-\alpha = - s_m^{-2} \sum_{k=1}^m \nu_{m,k}
\mathcal{E}_k - s_m^{-2} \sum_{k=1}^m \nu_{m,k} \log(L(x_k)) .
\end{gather}
Applying Lemma \ref{lem:logper}, we  obtain the following bound:
\begin{gather}\label{eq:msegph}
  \mathbb{E}\Big[\Big\{\sum_{k=1}^m \nu_{m,k}\log(\mathcal{E}_k)\Big\}^2\Big] \leq
  C(\delta,\eta^*) \,m.  
\end{gather}
The bias term is dealt with by applying Lemma \ref{lem:sumlll} wich yields
\begin{gather}\label{eq:biaisgph}
  \Big|\sum_{k=1}^m \nu_{m,k} \log(L(x_k))\Big| \leq m
  \eta^*(x_m)\{1+o(1)\},
\end{gather}
uniformly with respect to $|\eta|\leq \eta^*$.  Choosing $m$ as
in~(\ref{eq:condmgph}) yields~(\ref{eq:bornesup}). \qed

\section{Bandwidth selection} \label{sec:bandwidth}
In any semiparametric procedure, the main issue is the bandwidth
selection, here the number $m$ of Fourier frequencies used in the
regression. Many methods for choosing $m$ have been suggested, all
assuming some kind of second order regular variation of the spectral
density at~0. In Figures~\ref{fig:good}-~\ref{fig:horror} below, the
difficulty of choosing $m$ is illustrated, at least visually. In each
case, the values of the GPH estimator are plotted against the
bandwidth $m$, for values of $m$ between $10$ and $500$ and sample
size 1000.

In Figure~\ref{fig:good} the data is a simulated Gaussian
ARFIMA(0,$d$,0).  The spectral density $f$ of an ARFIMA(0,$d$,0)
process is defined by $f(x)=\sigma^2|1-\mathrm e^{\mathrm i
  x}|^{-2d}/(2\pi)$, where $\sigma^2$ is the innovation variance. Thus
it is second order regularly varying at zero and satisfies $f(x) =
x^{-\alpha} (C+O(x^2))$ with $\alpha=2d$. The optimal choice of the
bandwidth is of order $n^{4/5}$ and the semiparametric optimal rate of
convergence is $n^{2/5}$. Of course, it is a regular parametric model,
so a $\sqrt n$ consistent estimator is possible if the model is known,
but this is not the present framework. The data in
Figure~\ref{fig:medium} comes from an ARFIMA(0,$d$,0) observed in
additive Gaussian white noise with variance $\tau^2$. The spectral
density of the observation is then
\begin{gather*}
  \frac{\sigma^2}{2\pi} |1-\mathrm e^{\mathrm i x}|^{-2d} +
  \frac{\tau^2}{2\pi} = \frac{\sigma^2}{2\pi} |1-\mathrm e^{\mathrm i
    x}|^{-2d} \left\{ 1 + \frac{\tau^2}{\sigma^2} |1-\mathrm
    e^{\mathrm i x}|^{2d} \right\} \; .
\end{gather*}
It is thus second order regularly varying at 0 and the optimal rate of
convergence is $n^{2d/(4d+1)}$, with optimal bandwidth choice of order
$n^{4d/(4d+1)}$.  In Figures~\ref{fig:good} and~\ref{fig:medium} the
outer lines are the 95\% confidence interval based on the central
limit theorem for the GPH estimator of $d=\alpha/2$. See
\cite{robinson:1995l}.

A visual inspection of Figure~\ref{fig:good} leaves little doubt that
the true value of $d$ is close to .4. In Figure~\ref{fig:medium}, it
is harder to see that the correct range for the bandwidth is somewhere
betwen 50 and 100. As it appears here, the estimator is always
negatively biased for large $m$, and this may lead to underestimating
the value of $d$. Methods to correct this bias (when the correct
model is known) have been proposed and investigated
by~\cite{hurvich:ray:2003} and~\cite{hurvich:moulines:soulier:2005},
but again this is not the framework considered here.

Finally, in Figure~\ref{fig:horror}, the GPH estimator is computed for
a Gaussian process with autocovariance function $\gamma(k) = 1/(k+1)$
and spectral density $\log|1-\mathrm e^{\mathrm ix}|^2$. The true
value of $\alpha$ is zero, but the spectral density is infinite at
zero and slowly varying. The plot $\hat d(m)$ is completely
misleading. This picture is similar to what is called the Hill
``horror plot'' in tail index estimation. The confidence bounds are
not drawn here because there are meaningless. See
Example~\ref{xmpl:convproba}.

There has recently been a very important literature on methods to improve the
rate of convergence and/or the bias of estimators of the long memory parameter,
always under the assumption of second order regular variation. If this
assumption fails, all these methods will be incorrect. It is not clear if it is
possible to find a realistic method to choose the bandwidth $m$ that would
still be valid without second order regular variation. It might be of interest
to investigate a test of second order regular variation of the spectral
density.

\begin{figure}[htbp]
  \centering
  \includegraphics[width=12cm]{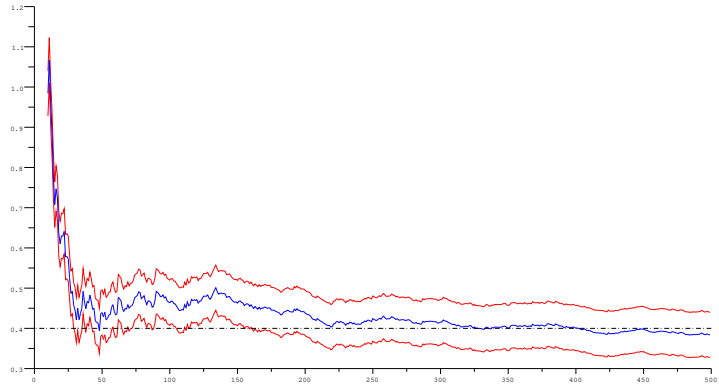}
  \caption{GPH estimator for ARFIMA(0,.4,0)}
  \label{fig:good}
\end{figure}

\begin{figure}[htbp]
  \centering
  \includegraphics[width=12cm]{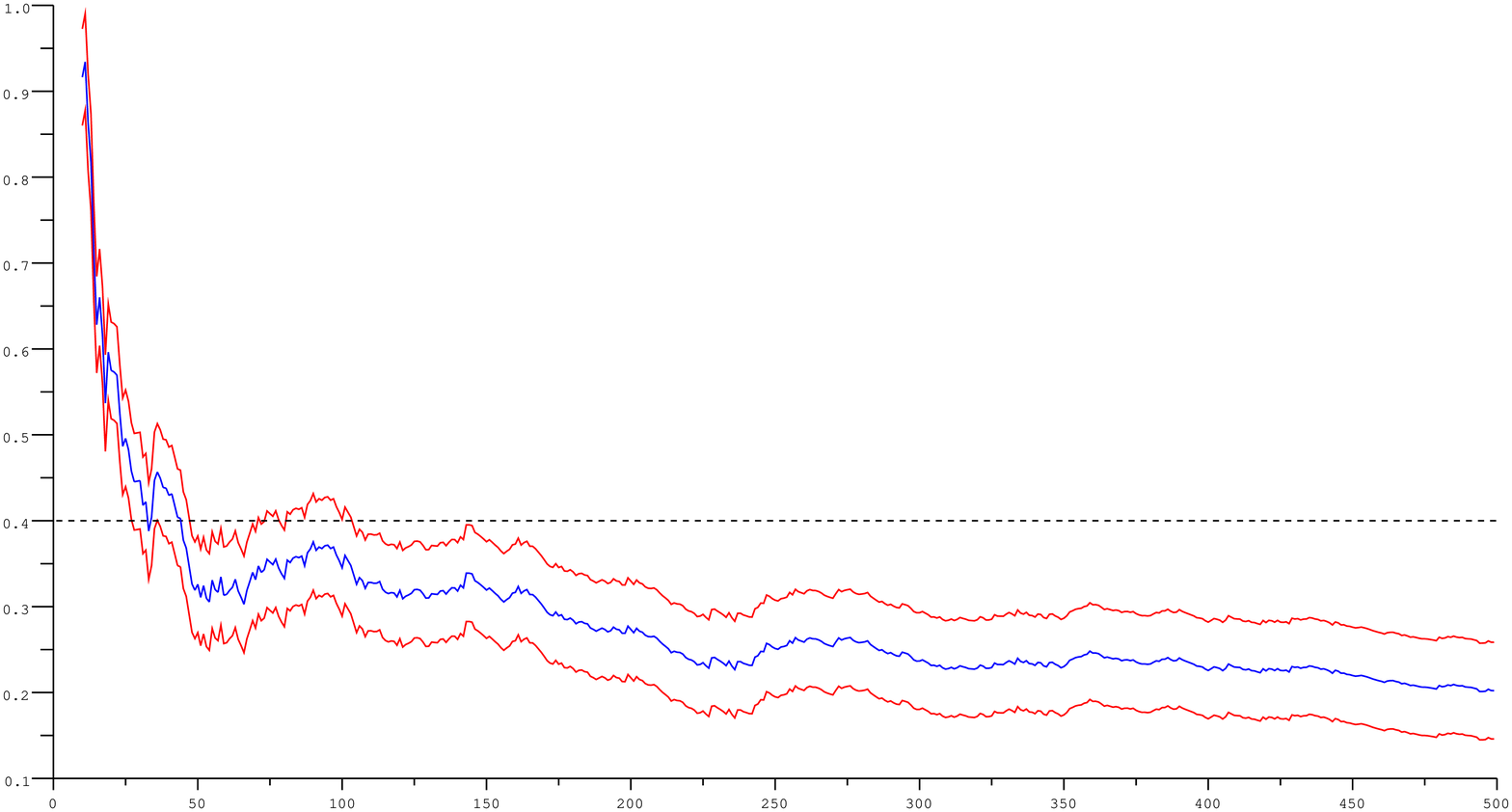}
  \caption{GPH estimator for ARFIMA(0,.4,0)+ noise}
  \label{fig:medium}
\end{figure}

\begin{figure}[htbp]
  \centering
  \includegraphics[width=12cm]{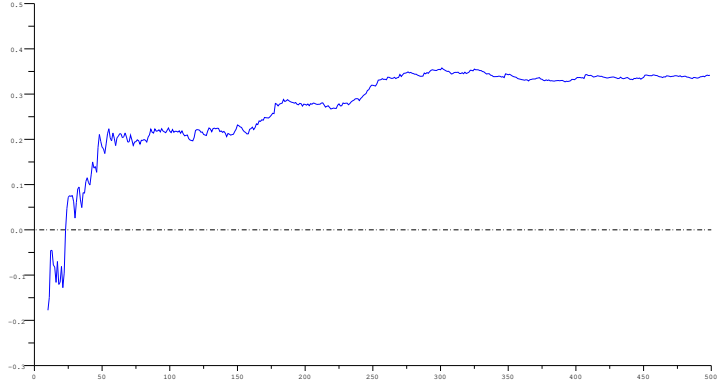}
  \caption{GPH ``horror'' plot}
  \label{fig:horror}
\end{figure}

\section{Technical  results} \label{sec:lemmes}

\begin{lemma} \label{lem:sumlll}
  Let $\eta^*$ be a non decreasing slowly varying function on
  $(0,\pi]$ such that $\lim_{s\to0}\eta^*(s)=0$.  Let $\eta$ be a
  measurable function on $(0,\pi]$ such that $|\eta|\leq \eta^*$, and
  define $h(x) = -\int_x^\pi \frac{\eta(s)}s \, \mathrm d s$ and
  $h^*(x) = \int_x^\pi \frac{\eta^*(s)}s \, \mathrm d s$.  Then, for
  any non decreasing sequence $m\leq n$,
 \begin{align}  
\label{eq:sumnuh}
& \left|\sum_{k=1}^m \nu_m(k) h(x_k)\right| \leq m 
   \eta^*(x_m) + O(\log^2(m)\eta^*(x_m) + \log(m) h^*(x_m)),
\end{align}
uniformly with respect to $|\eta|\leq\eta^*$.
\end{lemma}

\noindent{\em Proof.} \ 
Since $\eta$ is slowly varying, the function $h$ is also slowly
varying and satisfies $\lim_{x\to0} \eta(x)/h(x) = 0$. Then,
  \begin{align*}
    \sum_{k=1}^m h(x_k) & = \frac n {2\pi} \int_0^{x_m} h(s) \,ds + \frac n
    {2\pi} \sum_{k=1}^m \int_{x_{k-1}}^{x_k} \{h(x_k)-h(s)\} \,ds \\
    & = \frac n {2\pi} \int_0^{x_m} h(s) \,ds + \frac n {2\pi} \sum_{k=1}^m
    \int_{x_{k-1}}^{x_k} \int_s^{x_k} \frac{\eta(t)}t \, dt \,ds \; .
  \end{align*}
  Thus, for $|\eta| \leq \eta^*$ and $\eta^*$ increasing, 
  \begin{align*}
    \left| \sum_{k=1}^m h(x_k) - \frac n {2\pi} \int_0^{x_m} h(s) \,ds
    \right| & \leq \eta^*(x_m) \frac n {2\pi} \sum_{k=1}^m
    \int_{x_{k-1}}^{x_k} \int_s^{x_k} \frac{dt}t  \ ds \\
    & = \eta^*(x_m)
    \left\{ \sum_{k=1}^m \log(x_k) -  \frac n {2\pi} \int_0^{x_m} \log(s) \, ds\right\} \\
    & = \eta^*(x_m)
    \left\{ \sum_{k=1}^m \log(k) - m \log(m)+m\right\} \\
    & = O(\eta^*(x_m)\log(m)) \; .
  \end{align*}
  By definition, it holds that $xh'(x) = \eta(x)$.  Integration by
  parts yield
 \begin{gather} \label{eq:ipp}
   \frac n {2\pi} \int_0^{x_m} h(s) \,ds = mh(x_m) - \frac n {2\pi}
   \int_0^{x_m} \eta(s) \, \mathrm d s \; .
 \end{gather}
 Thus, 
\begin{gather}
  \label{eq:borneinth}
  \sum_{k=1}^m h(x_k) = m h(x_m) - \frac n {2\pi} \int_0^{x_m} \eta(s)
  \, \mathrm d s + O(\log(m)\eta^*(x_m)) \; .
\end{gather}
Similarly, we have:
  \begin{multline}
    \sum_{k=1}^m h(x_k) \log(x) = m h(x_m)\log(x_m) - \frac n {2\pi}
    \int_0^{x_m} \{\eta(s) \log(s) + h(s)\} \, \mathrm d x \\
    + O(\log(m)\{\log(x_m)\eta^*(x_m)+h^*(x_m)\}) \; .
    \label{eq:expressionsumhlog}
  \end{multline}
By definition of $\nu_m(k)$, we have:
\begin{align}\nonumber
  \nu_m(k) & = \log(k) - \frac1m\sum_{j=1}^m \log(j) = \log(x_k) -
  \log(x_m) + 1 + O\left(\frac{\log(m)}m\right).
\end{align}
Hence, applying \eqref{eq:ipp}, \eqref{eq:borneinth} and
\eqref{eq:expressionsumhlog}, we obtain
\begin{align*}
  \sum_{k=1}^m & \nu_m(k) h(x_k) &
  \\
  & = \sum_{k=1}^m \log(x_k) h(x_k) - \{\log(x_m) - 1 +
  O(\log(m)/m)\}\sum_{k=1}^m h(x_k)
  \\
  & = \frac n{2\pi} \int_0^{x_m} \eta(s) \{\log(x_m)-\log(s)\} \,
  \mathrm d s
  \\
  & \ \ \ \ \ \ \ \ \ \ \ \ \ \ \ \ \ \ \ \ \ \ \ \ \ \ \ \ \ \ +
  O(\log(m)\{\eta^*(x_m) \log(x_m) + h^*(x_m)\}) \; .
\end{align*}
Finally, since $|\eta|\leq \eta^*$ and $\eta^*$ is non decreasing, we
obtain
\begin{multline*}
  \frac n{2\pi} \left| \int_0^{x_m} \eta(s) \{\log(x_m) -
    \log(s)\}\,ds\right| \\
  \leq \frac n{2\pi} \eta^*(x_m) \int_0^{x_m} \int_s^{x_m} \frac 1 t
  \, dt \,ds = m\eta^*(x_m).
\end{multline*}
This yields \eqref{eq:sumnuh}.  \qed

\begin{lemma}\label{lem:logper}
  Let $\eta^*$ be a non decreasing slowly varying function such that
  $\lim_{x\to0}\eta^*(x)=0$.  Let $X$ be a {\em Gaussian} process with
  spectral density $f(x)=x^{-\alpha}L(x)$, where
  $\alpha\in[-\delta,\delta]$ and $L\in SV(\eta^*)$.  Let
  $\gamma=0,577216...$ denote Euler's constant.  Then, for all $n$ and
  all $k,j$ such that $0 < x_k, \, x_j \leq \pi/2$,
\begin{gather*}
  \left|\mathbb{E}[\log(\mathcal{E}_k)]-\gamma\right| + \left|\mathbb{E}[\log^2(\mathcal{E}_k)]
    - \frac{\pi^2}6\right| \leq C(\delta,\eta^*)\log(1+k) k^{-1}, \\
  |\mathbb{E}[\log(\mathcal{E}_k)-\gamma)(\log(\mathcal{E}_j)-\gamma)]| \leq
  C(\delta,\eta^*) \log^2(j) k^{-2}.
\end{gather*}
\end{lemma}

\paragraph{\bf Proof of Lemma \ref{lem:logper}}
It is well known (see for instance \cite{hurvich:deo:brodsky:1998},
\cite{moulines:soulier:1999}, \cite{soulier:2001}) that the bounds of
Lemma~\ref{lem:logper} are consequences of the covariance inequality
for functions of Gaussian vectors of \cite[Lemma~1]{arcones:1994} and
of the following bounds. For all $n$ and all $k,j$ such that $0 <
|x_k| \leq |x_j| \leq \pi/2$,
\begin{multline*}
  |\mathrm{cov}(d_{X,k}, d_{X,j})| + |\mathrm{cov}(d_{X,k}, \bar d_{X,j})-f(x_k)
  \delta_{k,j}|
  \\
  \leq C(\delta,\eta^*) \sqrt{f(x_k)f(x_j)}\, \log(j)\, k^{-1} \; .
\end{multline*}
Such bounds have been obained when the spectral density is second
order regularly varying. We prove these bounds under our assumptions
that do not imply second order regular varition. Denote $D_n(x) = (2
\pi n)^{-1/2} \sum_{t=1}^n \mathrm{e}^{-\mathrm{i}tx}$. Then
\[
\mathrm{cov}(d_{X,k},d_{X,j}) = \int_{-\pi}^\pi f(x) D_n(x_k-x)D_n(x_j+x) dx.
\]
Recall that by definition of the class $SV(\eta^*)$, there exists a
function $\eta$ such that $|\eta|\leq \eta^*$ and $L(x) = L(\pi)
\exp\{-\int_x^\pi s^{-1}\eta(s) \mathrm d s\}$. Since only ratio
$L(x)/L(\pi)$ are involved in the bounds, without loss of generality,
we can assume that $L(\pi)=1$.
We first prove that for all $k$ such that $x_k\leq \pi/2$,
\begin{gather}   \label{eq:keyresult1} 
  \left | \int_{-\pi}^\pi \left(\frac{f(x)}{f(x_k)} -1 \right)|
    D_n(x_k-x)|^2 \,\mathrm d x \right | \leq
  C(\delta,\eta^*) \log(k) k^{-1} \; .
\end{gather}
Since $L \in SV(\eta^*)$, the functions $x^\epsilon L(x)$ and
$x^\epsilon L^{-1}(x)$ are bounded for any $\epsilon>0$ and
\begin{gather}\label{eq:borneh} 
  \sup_{x\in[0,\pi]} x^{\delta}(L(x)+ L^{-1}(x)) \leq
  C(\eta^*,\delta) \; , \\
  \sup_{\alpha\in[-1+\delta,1-\delta]} \int_{-\pi}^\pi f(x) \, \mathrm
  d x \leq C(\eta^*,\delta) \ . \label{eq:borneintphi}
\end{gather}
Since $\eta^*$ is increasing, for all $0 < x < y \leq \pi/2$, it
holds that
\begin{align*}
  |f(x) - f(y)| & = |x^{-\alpha} L(x) - y^{-\alpha} L(y)| \\
  &  \leq \int_x^y |\alpha -  \eta(s)| s^{-\alpha-1} L(s) \,\mathrm d s \\
  & \leq \int_x^y (1 + \eta^*(\pi)) s^{-\alpha-1} L(s) \,\mathrm d
  s \; .
\end{align*}  
Since $\alpha\in[-1+\delta,1-\delta]$, $x^{-\alpha-1}L(x)$ is
decreasing.  Hence
\begin{gather}\label{eq:lipschitz}
  |f(x) - f(y)| \leq C(\eta^*,\delta) x^{-1} f(x)(y-x) \; .
  \end{gather}
  Define $F_n(x)= |D_n(x)|^2$ (the Fejer kernel).  We have
\begin{align}\label{eq:fejer}
  & \sup_{\pi/2 \leq |x| \leq \pi} |F_n(x_k-x)| = O(n^{-1}) \; , \ \ 
  \ \ \int_{-\pi}^{\pi} F_n(x) \, \mathrm d x = 1 \; , \\
  \label{eq:convoldirichlet}
  & \int_{-\pi}^{\pi}D_n(y+x)D_n(z-x)\, \mathrm d x=(2\pi n)^{-1/2}D_n(y+z) \; .
\end{align}  
From now on, $C$ will denote a generic constant which depends only on
$\eta^*$, $\delta$ and numerical constants, and whose value may
change upon each appearance. Applying \eqref{eq:borneh},
\eqref{eq:borneintphi} and \eqref{eq:fejer}, we obtain
\[
\int_{\pi/2 \leq |x| \leq \pi} \left| f^{-1}(x_k) f(x) - 1 \right|
F_n(x-x_k) \, \mathrm d x \leq C n^{-1} (f^{-1}(x_k) + 1) \leq C k^{-1} \; .
\]
The integral over $[-\pi/2,\pi/2]$ is split into integrals over
$[-\pi/2, -x_k/2]\cup[2x_k,\pi/2]$, $[-x_k/2,x_k/2]$ and $[x_k/2,2x_k]$.
If $x\in[-\pi/2, -x_k/2]\cup[2x_k,\pi/2]$, then $F_n(x-x_k)\leq C
n^{-1}x^{-2}$. Hence, applying Karamata's Theorem (cf.
\cite[Theorem~1.5.8]{bingham:goldie:teugels:1989}), we obtain:
\begin{align*}
  \int_{-\pi/2}^{-x_k/2} & + \int_{-\pi/2}^{-x_k/2} f(x) F_n(x) \, \mathrm d x
  \\
  & \leq C n^{-1}\int_{x_k/2}^{\pi/2} x^{-\alpha-2} L(x) \, \mathrm d x \leq
  Cn^{-1} x_k^{-\alpha-1} L(x_k) \leq C k^{-1} f(x_k) \; ,
  \\
  \int_{-\pi/2}^{-x_k/2} & + \int_{-\pi/2}^{-x_k/2} F_n(x) \, \mathrm
  d x \leq C n^{-1}\int_{x_k/2}^\infty x^{-2} \, \mathrm d x \leq Cn^{-1}
  x_k^{-1} \leq C k^{-1} \; .
\end{align*}
For $x \in [-x_k/2,x_k/2]$, $F_n(x_k-x) \leq n^{-1}x_k^{-2}$.  Thus,
applying again Karamata's Theorem, we obtain:
\begin{align*}
  \int_{-x_k/2}^{x_k/2} f(x) F_n(x-x_k) \, \mathrm d x & \leq C n^{-1} x_k^{-2}
  \int_{-x_k/2}^{x_k/2} x^{-\alpha} L(x) \, \mathrm d x
  \\
  & \leq C n^{-1} x_k^{-2} x_k^{-\alpha+1} L(x_k) \leq C k^{-1} f(x_k)
  \; ,
  \\
  \int_{-x_k/2}^{x_k/2} F_n(x-x_k)dx & \leq C n^{-1} x_k^{-1} \leq C
  k^{-1} \; .
\end{align*}
Applying \eqref{eq:lipschitz} and the bound $\int_{-x_k/2}^{x_k}
|x|F_n(x)| \, \mathrm d x \leq C n^{-1} \log(k)$, we obtain:
\begin{multline*}
  \int_{x_k/2}^{2x_k} |f(x) -f(x_k)|  F_n(x-x_k) \, \mathrm d x \\
  \leq C x_k^{-\alpha-1} L(x_k/2) \int_{x_k/2}^{2x_k}
  |x-x_k|F_n(x-x_k)| \, \mathrm d x \leq C f(x_k) k^{-1}\log(k) \; .
\end{multline*}
This proves \eqref{eq:keyresult1}.
We now prove that  all $k, j$ such that $0 < x_k \ne |x_j| \leq \pi/2$,
\begin{multline}\label{eq:keyresult2} 
  \left | \int_{-\pi}^{\pi} \left( \frac{f(x)}{f(x_k)} - 1
    \right) D_n(x_k -x) \overline{D_n(x_j-x)} \, \mathrm d x \right | \\
  + \left | \int_{-\pi}^{\pi} \left( \frac{f(x)}{f(x_k)}
      - 1 \right) D_n(x_k -x) D_n(x_j-x) \, \mathrm d x \right |\\
  \leq C(\delta,\eta^*) \log(k \vee |j|) (k \wedge |j|)^{-1}.
\end{multline}
Define $E_{n,k,j}(x) := D_n(x_k -x)
\overline{D_n(x_j-x)}$.  Since $0 \leq x_k, x_j \leq \pi/2$, for
$\pi/2 \leq |x| \leq \pi$, we have $|E_{n,k,j}(x)| \leq C n^{-1}$.
Hence, as above,
\[
\int_{\pi/2 \leq |x| \leq \pi} | f^{-1}(x_k)f(x)-1| \, \mathrm d x \leq C n^{-1}
(x_k^{\alpha} L^{-1}(x_k) + 1) \leq C k^{-1}.
\]
We first consider the case $k<j$ and we split the integral over
$[-\pi/2,\pi/2]$ into integrals over
$[-\pi/2,-x_k/2]\cup[2x_j,\pi/2]$, $[-x_k/2,x_k/2]$,
$[x_k/2,(x_k+x_j)/2]$, $[(x_k+x_j)/2, 2x_j]$, 
denoted respectively $I_1$, $I_2$, $I_3$ and $I_4$.\\
$\bullet$ The bound for the integral over
$[-\pi/2,-x_k/2]\cup[2x_j,\pi/2]$ is obtained as above (in the case
$k=j$) since $|E_{n,k,j}|\leq C n^{-1} x^{-2}$.  Hence $|I_1| \leq C
k^{-1}$.\\
$\bullet$ For $x \in [-x_k/2,x_k/2]$, $|E_{n,k,j}(x)| \leq C n^{-1} x_k^{-2}$,
hence we get the same bound: $|I_2| \leq C k^{-1}$.\\
$\bullet$ To bound $I_3$, we note  that on the interval $[x_k/2,(x_k+x_j)/2]$,
\[
|E_{n,k,j}(x)| \leq Cn^{1/2}(j-k)|D_n(x-x_k)|,
\]
and $ n^{1/2} |x-x_k| |D_n(x-x_k|$ is uniformly bounded.
Hence, applying \eqref{eq:lipschitz}, we obtain 
\[|I_3| \leq C(j-k)^{-1} x_k^{-1}x_j \leq C k^{-1}.\]
$\bullet$ The bound for $I_4$ is obtained similarly: $|I_4|\leq C k^{-1}\log(j)$.

\noindent $\bullet$ To obtain the bound in the case $x_j<x_k$, the interval $[-\pi,\pi]$
is split into $[-\pi,-\pi/2]\cup[\pi/2\pi]$ $[-\pi/2,-x_k/2]\cup[2x_k,
\pi/2]$, $[-x_k/2,x_j/2]$, $[x_j/2,(x_k+x_j)/2]$ and $[(x_k+x_j)/2,
2x_k]$.  The arguments are the same except on the interval
$[-x_k/2,x_j/2]$ where a slight modification of the argument is
necessary.  On this interval, it still holds that $|E_{n,k,j}(x)|\leq
n^{-1}x_k^{-2}$.  Moreover, $x^{\delta} L(x)$ can be assumed
increasing on $[0,x_k/2]$, and we obtain:
\[
\int_{-x_k/2}^{x_j/2} x^{-\alpha} L(x) \, \mathrm d x \leq x_k^{\delta} L(x_k)
\int_{-x_k/2}^{x_j/2} x^{-\alpha-\delta} \, \mathrm d x \leq C x_k^{-\alpha+1}
L(x_k) \; .
\]
The rest of the argument remains unchanged.

\noindent $\bullet$ To obtain the bound in the case $ x_j < 0 < x_k$, the interval
$[-\pi,\pi]$ is split into $[-\pi,-\pi/2]\cup[\pi/2,\pi]$
$[-\pi/2,2x_j]\cup[2x_k, \pi/2]$, $[2x_j,-x_k/2]$, $[-x_k/2,x_k/2]$
and $[x_k/2, 2x_k]$ and the same arguments are applied.
\qed


\end{document}